\documentclass[12pt]{amsart}
\usepackage{amsmath,amsthm,amsfonts,amssymb}
\begin{document} 
\newcommand{\B}{{\mathbb B}}
\newcommand{\C}{{\mathbb C}}
\newcommand{\Ok}{{\mathcal O}_k}
\renewcommand{\O}{{\mathcal O}}
\newcommand{\Ol}{{\mathcal O}_L}
\newcommand{\N}{{\mathbb N}}
\newcommand{\F}{{\mathbb F}}
\newcommand{\Q}{{\mathbb Q}}
\renewcommand{\H}{{\mathbb H}}
\newcommand{\Z}{{\mathbb Z}}
\renewcommand{\P}{{\mathbb P}}
\newcommand{\R}{{\mathbb R}}
\newcommand{\rc}{\subset}
\newcommand{\rank}{\mathop{rank}}
\newcommand{\Tr}{\mathop{Tr}}
\newcommand{\dimc}{\mathop{dim}_{\C}}
\newcommand{\tensor}{\otimes}
\newcommand{\Lie}{\mathop{Lie}}
\newcommand{\Auto}{\mathop{{\rm Aut}_{\mathcal O}}}
\newcommand{\alg}[1]{{\mathbf #1}}
\newtheorem*{definition}{Definition}
\newtheorem*{claim}{Claim}
\newtheorem{corollary}{Corollary}
\newtheorem{conjecture}{Conjecture}
\newtheorem*{Sconjecture}{Schanuel's Conjecture}
\newtheorem*{SpecAss}{Special Assumptions}
\newtheorem{example}{Example}
\newtheorem*{remark}{Remark}
\newtheorem*{observation}{Observation}
\newtheorem*{fact}{Fact}
\newtheorem*{remarks}{Remarks}
\newtheorem{lemma}{Lemma}
\newtheorem{proposition}{Proposition}
\newtheorem{theorem}{Theorem}
\title[Elliptic curves in $SL_2(\C)/\Gamma$]{%
On Elliptic Curves in $SL_2(\C)/\Gamma$, Schanuel's conjecture
and geodesic lengths
}
%
\author {J\"org Winkelmann}
\begin{abstract}
Let $\Gamma$ be a discrete cocompact subgroup of $SL_2(\C)$.
We conjecture that the quotient manifold $X=SL_2(\C)/\Gamma$
contains infinitely many non-isogenous elliptic curves
and prove this is indeed the case if Schanuel's conjecture
holds. We also prove it in the special case where
$\Gamma\cap SL_2(\R)$ is cocompact in $SL_2(\R)$.

Furthermore, we deduce some consequences for the geodesic
length spectra of real hyperbolic $2$- and $3$-folds.
\end{abstract}
\address{%
 J\"org Winkelmann \\
 Institut Elie Cartan (Math\'ematiques)\\
 Universit\'e Henri Poincar\'e Nancy 1\\
 B.P. 239, \\
 F-54506 Vand\oe uvre-les-Nancy Cedex,\\
 France
}
\email{jwinkel@member.ams.org\newline\indent{\itshape Webpage: }%
http://www.math.unibas.ch/\~{ }winkel/
}
\thanks{
{\em Acknowledgement.}
The author wants to thank the University of Tokyo
and the Korea Institute for Advanced Study in Seoul.
The research for this article was done during the stays of the
author at these institutes.}
\maketitle
\section{Introduction}
Let $\Gamma$ be a discrete cocompact subgroup of $SL_2(\C)$.
We are interested in closed complex analytic subspaces of the 
complex quotient manifold $X=SL_2(\C)/\Gamma$.
It is well-known that $X$ contains no hypersurfaces and it is
easy to show that it contains no curves of genus $0$.
The existence of curves of genus $\ge 2$ is an unsolved problem.

On the other hand, 
it is not hard to show that there do exist curves of genus one
(elliptic curves).
(For these assertions, see \cite{HW},\cite{SMF}.)

Our goal is to investigate how many different
curves of genus one can be embedded in one such quotient manifold.
There are only countably many abelian varieties which can be embedded
into a quotient manifold of a complex semisimple Lie group by a discrete
cocompact subgroup (\cite{SMF},~Cor.~4.6.2).
Thus the question is: Is the number of non-isomorphic elliptic
curves in such a quotient $SL_2(\C)/\Gamma$ finite or countably infinite?

Under the additional assumption that $\Gamma\cap SL_2(\R)$ is cocompact
in $SL_2(\R)$ we show that there are infinitely many isogeny classes
of elliptic curves in $X$ (thm.~\ref{prop-real}).
We will see that there do exist discrete cocompact subgroups in
$SL_2(\C)$ with this property (cor.~\ref{cor-ex-real}).
We conjecture that this additional assumption
($\Gamma\cap SL_2(\R)$ being cocompact
in $SL_2(\R)$) is not needed and show that it can be dropped provided
{\em Schanuel's conjecture} is true (see cor.~\ref{cor-schanuel}).

In order to show that there are infinitely many non-isogenous
elliptic curves, one first has to construct elliptic curves
and then one has to investigate under which conditions they are
isogenous.
There is a well-known way to construct elliptic curves in $X=SL_2(\C)/\Gamma$, 
going back to ideas of Mostow (\cite{Mos}). In fact
every elliptic curve in $X$ arises in this way (\cite{HW}).
This method works as follows:
If $\gamma\in\Gamma$ is a semisimple element of infinite order,
then the centralizer $C=\{g\in SL_2(\C):g\gamma=\gamma g\}$
is isomorphic to $\C^*$ as a complex Lie group and $C\cap\Gamma$
is a discrete subgroup containing $\gamma$ and therefore
commensurable with $\{\gamma^k:k\in\Z\}$.
The quotient of $\C^*$ by an infinite discrete subgroup is necessarily
compact. Hence for every semisimple element $\gamma\in\Gamma$ 
of infinite order we obtain an elliptic curve $E\subset X=SL_2(\C)/\Gamma$
which arises as orbit of the centralizer $C$. Moreover,
this elliptic curve $E\simeq C/(C\cap\Gamma)$ is isogenous to
$C/\left<\gamma\right>$ and therefore isogenous to
$\C^*/\left<\lambda\right>$ where $\lambda$ and $\lambda^{-1}$
are the eigenvalues of the matrix  $\gamma\in SL_2(\C)$.

Thus our problem is to investigate how many different eigenvalues
occur and under which circumstances different eigenvalues
lead to non-isogenous elliptic curves.

First we show that for every Zariski-dense subgroup
$\Gamma\subset SL_2(\C)$ there are infinitely many pairwise
multiplicatively independent complex numbers occuring as eigenvalues
for elements of $\Gamma$ (thm.~\ref{thm-mult-ind}).

We conjecture that, if the eigenvalues are algebraic numbers
(this is known to be the case if $\Gamma$ is cocompact), then
multiplicatively independent eigenvalues always lead to
non-isogenous elliptic curves. We can prove that this conjecture
holds if Schanuel's conjecture from transcendental number theory
is true.

Even without assuming Schanuel's conjecture to be true we can
prove the existence of infinitely many non-isogenous
elliptic curves in the case where the eigenvalues are real.

In this way we obtained the desired result in the special
case where the intersection $\Gamma\cap SL_2(\R)$ is
cocompact in $SL_2(\R)$.

Using an arithmetic construction one can show that
discrete cocompact subgroups $\Gamma$ for which $\Gamma\cap SL_2(\R)$
is cocompact in $SL_2(\R)$ do indeed exist.

These results on elliptic curves in $SL_2(\C)/\Gamma$ can be related
to questions on the length of closed geodesics on real hyperbolic
manifolds of dimension $2$ or $3$.
More precisely, let $M$ be a compact real Riemannian manifold (without
boundary) of dimension $2$ or $3$ which carries a Riemannian
metric of constant negative curvature.
Let $\Lambda$ be set of all positive real numbers occuring as
length of a closed geodesic on $M$. Then $\Lambda$ contains infinitely
many elements which are pairwise linearly independent over $\Q$
(thm.~\ref{thm-geodesic}).
\section{Multiplicatively independent eigenvalues}
\subsection{Announcement of theorem~\ref{thm-mult-ind}}
\begin{definition}\label{def-mult-dep}
Two non-zero elements $x,y$ in a field $k$ are called
{\em multiplicatively dependent}
if there exists a pair
$(p,q)\in\Z\times\Z\setminus\{(0,0)\}$ such that $x^q=y^p$.

They are called {\em multiplicatively independent}
if they are not multiplicatively dependent.
\end{definition}
By this definition a root of unity is multiplicatively dependent
with every other element of $k^*$.
Thus, if $x,y\in k^*$ are multiplicatively independent,
this implies in particular that neither $x$ nor $y$ is a root of unity.

Note that being multiplicatively dependent is an equivalence
relation on the set of all elements of $k^*$ which are not roots of unity.

The purpose to of this section is to prove the following
theorem:
\begin{theorem}\label{thm-mult-ind}
Let $\Gamma$ be a subgroup of $SL_2(\C)$ which is dense
in the algebraic Zariski topology.

Then there exists infinitely many pairwise multiplicatively
independent complex numbers $\lambda$
which occur as eigenvalues
for elements of $\Gamma$.
\end{theorem}

\subsection{A fact from Combinatorics}
As a preparation for the proof of thm.~\ref{thm-mult-ind}
we need a combinatorial fact.
\begin{lemma}\label{lemma-combin}
Let $S$ be a finite set, $\phi:\N\to S$ a map.

Then there exists a natural number $N\le\#S$ and an element $s\in S$
such that
\[
A_{s,N}=\{x\in\N: \phi(x)=s=\phi(x+N)\}
\]
is infinite.
\end{lemma}
\begin{proof}
Assume the contrary.
Then $A_{s,N}$ is a finite set for all $s\in S$,
$1\le N\le\#S$.
Hence there is a number $M\in\N$ such that $x<M$
for all $x\in\cup_{s\in S}\cup_{N\le\#S} A_{s,N}$.

But this implies that
$\phi(M+i)\ne\phi(M+j)$
for all $0\le i < j\le \#S$,
which is impossible by the pigeon-hole principle.
\end{proof}
\subsection{Roots in finitely generated fields}
We need the following well-known fact on finitely generated
fields.
\begin{lemma}\label{fin-gen-field}
Let $K$ be a finitely generated  field extension of $\Q$.

Then for every element $x\in K$ one of the properties
hold:
\begin{itemize}
\item
$x=0$,
\item
$x$ is an invertible algebraic integer (i.e.~a unit)
or
\item
there exists a discrete valuation $v:K^*\to \Z$ with $v(x)\ne 0$.
\end{itemize}

\end{lemma}
For the convenience of the reader we sketch a proof.
\begin{proof}
Let $K_0$ denote the algebraic closure of $\Q$ in $K$.
Then $K_0$ is a number field and $K$ can be regarded
as function field of a projective variety $V$ defined over $K_0$.
Let $f\in K$. If $f\not\in K_0$, then $f$ is a non-constant rational
function and therefore there is a discrete valuation given by
the pole/zero-order along a hypersurface which does not annihilate $f$.
If $f\in K_0$, then either $f=0$, or $f$ is a unit, 
i.e.~an invertible algebraic integer or
an extension of a $p$-adic valuation is non-zero for $f$.
\end{proof}

Let $K$ be a field and $W_K$ the group of roots of unity contained
in $K$. Let $x\in K^*$.
We want to measure up to which degree $d$ it is possible
to find a $d$-th root of $x$ in $K$ (modulo $W_K$).
For this purpose we 
define
\[
\rho_K(x)=\sup\{n\in\N:\exists \alpha\in K:\alpha^nx^{-1}\in W_K\}
\in\N\cup\{\infty\}.
\]
\begin{lemma}\label{lemma3}
Let $K$ be a finitely generated field extension of $\Q$ and $x\in K^*$. 
Then $\rho_K(x)<\infty$
unless $x$ is a root of unity.
\end{lemma}
\begin{proof}
Let $x$ be an element of $K^*$ which is not a root of unity.
First we discuss the case in which there exists a 
discrete valuation $v:K^*\to\Z$ with $v(x)\ne 0$.
In this case $\alpha^nx^{-1}\in W_K$ for $\alpha\in K$ implies
$v(\alpha)=\frac{1}{n}v(x)\in\Z$.
Therefore $\rho_K(x)\le|v(x)|$ in this case.

Now let us discuss the case where every discrete valuation 
on $K$ annihilates $x$.
By lemma~\ref{fin-gen-field}, this implies that $x$ is contained in the 
algebraic closure $K_0$ of $\Q$ in $K$ and moreover that $x\in\O_{K_0}^*$, 
i.e.{} $x$ is an invertible algebraic integer. 
Assume that there are elements $\alpha\in K$, $w\in W_K$ and
$n\in\N$ such that $\alpha^n=xw$.
Then $\alpha^{nN}=x^N$ for some $N\in\N$. As a consequence, $\alpha$
is integral over $\O_{K_0}$. Similarily, $\alpha^{-nN}=x^{-N}$
implies that $\alpha^{-1}$ is integral over $\O_{K_0}$.
Thus we obtain: {\em If $\alpha^nx^{-1}\in W_K$ for some $\alpha\in K$ and
$n\in\N$, then $\alpha\in\O^*_{K_0}$.}

Therefore
\[
\rho_K(x)=\sup\{n\in\N:\exists\alpha\in\O_{K_0}^*:
\alpha^nx^{-1}\in W_K\}.
\]
A theorem of Dirichlet states that $\O_{K_0}^*$ is a finitely
generated abelian group (with respect to multiplication).
Thus $\O_{K_0}^*/W_K\simeq\Z^d$ for some $d\in\N$.%
\footnote{More precisely, the theorem of Dirichlet states
$d=r+s-1$ where $r$ is the number of
real embeddings of $K_0$ and $s$ the number of pairs of conjugate
complex embeddings.}
This implies $\rho_K(x)<\infty$.
\end{proof}

\begin{lemma}\label{lemma5}
Let $K$ be a field, $x\in K^*$ with $\rho_K(x)<\infty$.
Assume that there are integers $p\in\Z$, $q\in\Z\setminus\{0\}$
and an element $\beta\in K^*$ such that $\beta^qx^{-p}\in W_K$.

Then $\frac{p}{q}\rho_K(x)\in\Z$.
\end{lemma}
\begin{proof}
Let $n=\rho_K(x)$. Assume that $\frac{p}{q}n\not\in\Z$ and
let $\Gamma$ denote the additive subgroup of $\Q$ generated by
$\frac{1}{n}$ and $\frac{p}{q}$. Now $\frac{1}{n}\Z\subsetneq\Gamma$, hence
there is  a natural number $N>n$ such that $\Gamma=\frac{1}{N}\Z$.
Since $\Gamma$ is generated by $1/n$ and $p/q$, there are integers
$k,m\in\Z$ such that
\[
k\frac{1}{n}+m\frac{p}{q}=\frac{kq+nmp}{nq}=\frac{1}{N}.
\]
Since $n=\rho_K(x)$, there is an element $\alpha\in K^*$ with
$\alpha^nx^{-1}\in W_K$.
Now we define
\[
\gamma=\alpha^k\beta^m.
\]
We claim that $\gamma^Nx^{-1}\in W_K$.
Indeed, since $\frac{1}{N}=nq/(kq+nmp)$, this condition is
equivalent to $\gamma^{nq}x^{-kq-nmp}\in W_K$
which can be verified as follows:
\[
\gamma^{nq}x^{-kq-nmp}=\alpha^{knq}\beta^{mnq}x^{-kq-nmp}
=\left(\alpha^nx^{-1}\right)^{kq}\left(\beta^qx^{-p}\right)^{nm}\in W_K.
\]
But $\gamma^Nx^{-1}\in W_K$ implies $\rho_K(x)\ge N$,
contradicting $N>n=\rho_K(x)$.

Thus we see that $\frac{p}{q}$ must be contained in $\frac{1}{n}\Z$.
\end{proof}

The statement of the lemma may be reformulated in the following way:
\begin{corollary}\label{cor-theta}
Let $K$ be a field, $x\in K^*$ with $\rho_K(x)<\infty$.
Let
\[
\Theta_{K,x}=\left\{\frac{p}{q}\in\Q:\exists\beta\in K^*: \beta^qx^{-p}\in W_K \right\}.
\]
Then $\Theta_{K,x}$ is a discrete subgroup of $(\Q,+)$, generated by
$\frac{1}{\rho_K(x)}$.
\end{corollary}

Next we verify that the behaviour of $\rho_K(x)$ under finite field
extensions is as to be expected.

\begin{lemma}\label{rho-finite}
Let $L/K$ be a finite field extension of degree $d$ and $x\in K^*$
with $\rho_K(x)<\infty$.

Then there exists a natural number $s$ which divides $d$ such that
$\rho_{L}(x)=s\rho_K(x)$.
\end{lemma}

\begin{proof}
In the notation of cor.~\ref{cor-theta} $\Theta_{K,x}$ is a subgroup
of $\Theta_{L,x}$.

On the other hand, if there is an element $\beta\in (L)^*$ and a
natural number
$n$ such that $\beta^nx^{-1}\in W_{L}$, then 
\[
N_{L/K}(\beta^nx^{-1})=
\left(N_{L/K}(\beta)\right)^nx^{-d}\in W_K
\]
 and consequently
$\frac{d}{n}\rho_K(x)\in\Z$ (lemma~\ref{lemma5}).
Thus $\frac{1}{n}\in\frac{1}{\rho_{L}(x)}\Z$ implies
$\frac{1}{n}\in\frac{1}{d\rho_K(x)}\Z$.

Combined, these facts yield 
\[
\frac{1}{\rho_K(x)}\Z 
\subset \frac{1}{\rho_{L}(x)}\Z
\subset \frac{1}{d\rho_K(x)}\Z.
\]
This implies the statement of the lemma.
\end{proof}
\subsection{An auxiliary proposition}
\begin{proposition}\label{lemma6}
Let $K$ be a finitely generated field extension of $\Q$, 
$\bar K$ an algebraic closure,
$S$ a finite subset of $K^*$
and $\Lambda\subset\bar K^*$ a subset
such that the following properties are fulfilled:
\begin{enumerate}
\item
$\deg K(\lambda)/K\le 2$ for every $\lambda\in\Lambda$,
\item
for every $\lambda\in\Lambda$ there exists an element $\mu\in S$
and integers $p,q\in\Z\setminus\{0\}$ such that
$\lambda^p=\mu^q$.
\end{enumerate}

Then there exists a finite subgroup $W\subset\bar K^*$ and
a finite subset $S'\subset\bar K^*$ such that
for every $\lambda\in\Lambda$ there exists an element
$\alpha\in S'$, an integer $N\in\Z$ and an element
$w\in W$ such that $\alpha^Nw=\lambda$.

Moreover, the set $S'$ can be chosen in such a way
that none of its elements is a root of unity.
\end{proposition}
\begin{proof}
For each element $\mu\in S$ which is not a root of unity
we choose an element $\alpha_\mu\in\bar K^*$
such that
\[
(\alpha_\mu)^{2\rho_K(\mu)}=\mu.
\]
Let $S'$ be the set of all these elements $\alpha_\mu$.
Evidently none of these elements $\alpha_\mu$ is a root of unity.
Let $L$ denote the field generated by $K$ and the elements of $S'$.
Note that $L$ is a finitely generated field. 
Let $L_0$ denote the algebraic closure of $\Q$ in $L$.
Then $L_0$ is a number field. Let $d_0$ denote its degree
(over $\Q$). Recall that for any natural number,
in particular for $2d_0$, there are only finitely
many roots of unity of degree not greater than this
number. Let $W$ 
be the set of all roots of unity $w$ in $\bar K^*$ for which
$\deg(L(w)/L)\le 2$. Then $\deg(w)\le 2d_0$
for every $w\in W$. Therefore $W$ is is a finite group. 
By construction it contains every root of unity which is in
$L(\lambda)$ for some $\lambda\in\Lambda$.

Now choose an arbitrary element $\lambda\in\Lambda$.
If $\lambda$ is a root of unity, it is contained in $W$ implying
that $\lambda=\alpha^0w$ for $w=\lambda$ and $\alpha$ arbitrary.
Thus we may assume that $\lambda$ is not a root of unity.
There are integers $p,q\in\Z\setminus\{0\}$ and an element
$\mu\in S$ such that $\lambda^p=\mu^q$.
Since $\lambda$ is not a root of unity, this implies that neither $\mu$
can be a root of unity. 
Thus $\rho_K(\mu)<\infty$ (lemma~\ref{lemma3}) 
and there is an element
$\alpha_\mu\in S'$ with $\alpha_\mu^{2\rho_K(\mu)}=\mu$.

By lemma~\ref{lemma5} the equality $\mu^q=\lambda^p$ implies
\[
\frac{q}{p}\rho_{K(\lambda)}(\mu)\in\Z.
\]
Thanks to lemma~\ref{rho-finite}
we know that either
$\rho_{K(\lambda)}(\mu)=\rho_K(\mu)$ or
$\rho_{K(\lambda)}(\mu)=2\rho_K(\mu)$.
In both cases it follows that
\[
2\frac{q}{p}\rho_K(\mu)\in\Z.
\]
In other words, there is an integer $N\in\Z$ such that
$2q\rho_K(\mu)=pN$.
Therefore
\[
(\alpha_\mu^N)^p=\alpha_\mu^{pN}=\alpha_\mu^{2q\rho_K(\mu)}=\mu^q=\lambda^p.
\]
Hence $\left(\alpha_\mu^{-N}\lambda\right)^p=1$.
Let $w=\alpha_\mu^{-N}\lambda$. Then $w$ is a root of unity
which is contained in the field $L(\lambda)$. It follows that $w\in W$.
Thus we have verified that there exist elements
$\alpha\in S'$, $N\in\Z$ and $w\in W$
such that $\alpha^Nw=\lambda$.
\end{proof}

\subsection{Proof of theorem \ref{thm-mult-ind}}
\begin{proof}
If $\Gamma$ is a Zariski-dense subgroup of $SL_2(\C)$,
then $\Gamma$ contains a finitely generated torsion-free
subgroup $\Gamma_0$
which is still Zariski-dense (see \cite{SMF}, lemma 1.7.12 and 
Prop.~1.7.2).
Fix a finite set $E$ of generators of $\Gamma_0$.
Let $k$ be the field generated by all the matrix coefficients of 
elements of $E$. Then $k$ is a finitely generated extension field
of $\Q$ and $\Gamma_0\subset SL_2(k)$.

Let $\Lambda$ be the set of all complex numbers other than $1$
and $-1$ occuring as an eigenvalue for an element $\gamma\in\Gamma_0$.
We observe that a number $\lambda\in\C^*\setminus\{1,-1\}$
is contained in $\Lambda$ if and only if there exists
an element $A\in\Gamma_0$ such that $\Tr(A)=\lambda+\lambda^{-1}$.
Since $\Gamma_0$ is Zariski dense, the set
\[
\left\{ \Tr(A): A\in\Gamma_0\right\}
\]
is Zariski dense in $\C$. 
It follows that $\Lambda$ is an infinite set.

We claim that $\Lambda$ contains no root of unity.
Indeed, assume that a root of unity $\omega$ is contained in $\Lambda$.
Then $\omega\ne 1,-1$ and consequently $\omega\ne\omega^{-1}$.
Therefore every element $A\in SL_2(\C)$ with $\omega$ as an eigenvalue
is conjugate to
\[
\begin{pmatrix} \omega & \\ & \omega^{-1}
                              \end{pmatrix}.
\]
As a consequence, such a matrix  $A$ is of finite order. This contradicts
the assumption that $\Gamma_0$ is torsion-free.
Thus $\Lambda$ can not contain any root of unity.

Let $\Sigma$ denote the set of all complex numbers which are roots of unity.
As remarked before,
the notion of ``multiplicative dependence'' defines an equivalence relation
on $\C^*\setminus\Sigma$.

Let us assume that the statement of the theorem fails.
Since $\Lambda\subset\C^*\setminus \Sigma$ and since ``multiplicative dependence''
defines an equivalence relation on $\C^*\setminus\Sigma$,
it follows that there is a finite set $S$ and complex numbers 
$(\mu_i)_{i\in S}\in \C^*\setminus\Sigma$ such that for every 
$\lambda\in\Lambda$
there exists an index $i\in S$ and non-zero integers
$p,q\in\Z\setminus\{0\}$ with $\lambda^p=\mu_i^q$.

Let $K$ denote the field generated by $k$ and all the elements
$\mu_i$ ($i\in S)$.
Recall that every element of $\Lambda$ is an eigen value
for a matrix in $SL_2(k)\subset SL_2(K)$.
Therefore $\deg (K(\lambda)/K)\le 2$ for every $\lambda\in\Lambda$.

We may now invoke proposition~\ref{lemma6}.

Thus we obtain  the following statement:
{\sl
There are finitely many complex numbers $(\alpha_i)_{i\in S}$,
none of which is a root of unity,
and a finite subgroup $W$ of the multiplicative group $\C^*$
such that
for every $\lambda\in\Lambda$ there are $i\in S$, $n\in\Z$
and $w\in W$ such that $\lambda=\alpha_i^nw$.}

By adjoining all the elements of $W$ to $K$, we also may deduce
that in this case there exists a finitely generated field $L$ containing all
the $\alpha_i$ ($i\in S$) and all
$\lambda\in\Lambda$ and $w\in W$.

Let $\lambda\in\Lambda$, $\zeta\in S$, $q\in\Z\setminus\{0\}$
and $w_0\in W$ such that 
$\lambda=w_0\alpha_\zeta^q$.

Then, after
replacing $\Gamma_0$ by $g\Gamma g^{-1}$ for an appropriately
chosen $g\in SL_2(\C)$,
we obtain
\[
\Gamma_0 \ni \gamma=
\begin{pmatrix}
\lambda & \\ & \lambda^{-1} 
\end{pmatrix}
=
\begin{pmatrix}
w_0\alpha_\zeta^q & \\ & w_0^{-1}\alpha_\zeta^{-q}
\end{pmatrix}.
\]

By the assumption of Zariski density  $\Gamma_0$ must also
contain an
element $\delta\in\Gamma_0$
which does not commute with $\gamma$. 
Let
\[
\delta=
\begin{pmatrix}
a & c \\ b & d 
\end{pmatrix}
\] be such an element.
By the assumption of Zariski density of $\Gamma_0$ we may and do
require that $a,d\ne 0$.

Let $g_n=\gamma^n\delta$ for $n\in\N$.

Using lemma~\ref{lemma-combin}, we conclude that there exists a natural number $N$,
 an infinite subset
$A\subset\N$, an index $\xi$, an element $\tilde w\in W$ 
and sequences of non-zero integers $m_k,m_k'\in\Z\setminus\{0\}$
such that
$\tilde w\alpha_\xi^{m_k}$ resp. $\tilde w\alpha_\xi^{m_k'}$ is an eigenvalue
of $g_k$ resp. $g_{k+N}$ for all $k\in A$.
Moreover, we may assume that all the numbers $m_k$ and $m'_k$ 
have the same sign.

Since $w_0$ is a root of unity,
we may invoke the pigeon-hole principle 
in order to deduce that (by replacing $A$ with an appropriate
smaller set) we may 
assume that there is an element $w_1\in W$ such that $w_0^k=w_1$ for
all $k\in A$. Let $w_2=w_1w_0^N$. Then $w_2=w_0^{k+N}$
for all $k\in\N$.

Now recall that for an element $g\in SL_2(\C)$ with eigenvalues
$\lambda,\lambda^{-1}$ we have $Tr(g)=\lambda+\lambda^{-1}$.

It follows that
\begin{equation}\label{eq-tr-1}
Tr(\gamma^k\delta)=
w_1\alpha_\zeta^{qk} a + w_1^{-1}\alpha_\zeta^{-qk}d 
= \tilde w\alpha_\xi^{m_k} + \tilde w^{-1}\alpha_\xi^{-m_k}
\end{equation}

and
\begin{equation}\label{eq-tr-2}
Tr(\gamma^{k+N}\delta)
=w_2\alpha_\zeta^{q(k+N)} a + w_2^{-1}\alpha_\zeta^{-q(k+N)}d 
= \tilde w\alpha_\xi^{m_k'} + \tilde w^{-1}\alpha_\xi^{-m_k'}
\end{equation}
for all $k\in A$.

Recall that $\alpha_\zeta$ is contained in the finitely generated field
$L$ and is not a root of unity.
Therefore there exists an absolute value
$|\ |$ on $L$ such that $|\alpha_\zeta|\ne 1$.
In what follows,
$|\ |$ always denotes this (possibly non-archimedean)
absolute value on $L$.

Using $|\alpha_\zeta|\ne 1$ and $a,d,q\ne 0$ we obtain
\[
\lim_{k\to\infty}\left| w_1\alpha_\zeta^{qk} a + w_1^{-1}\alpha_\zeta^{-qk}d 
\right | = +\infty
\]
Combined with eq.~(\ref{eq-tr-1}), this yields
\[
\lim_{k\to\infty}\left| 
 \tilde w\alpha_\xi^{m_k} + \tilde w^{-1}\alpha_\xi^{-m_k}
\right | = +\infty
\]
This is only possible if $|\alpha_\xi|\ne 1$.

Without loss of generality we may assume that $|\alpha_\zeta|,|\alpha_\xi|>1$,
$q>0$ and $m_k,m_k'>0$ for all $k\in A$.

Then
\[
\lim_{k\to\infty}\alpha_\zeta^{-qk}=0=\lim_{k\to\infty}\alpha_\xi^{-m_k}
                                     =\lim_{k\to\infty}\alpha_\xi^{-m'_k}.
\]
It follows that  the quotient of the respective left hand sides of the
equations (\ref{eq-tr-2}) and (\ref{eq-tr-1}) converges to
$
\alpha_\zeta^{qN}\frac{w_2}{w_1}$.
Evidently the quotient of the respective right hand
sides converges to the same value. Hence:
\[
\alpha_\zeta^{qN}\frac{w_2}{w_1}=\lim_{k\to\infty,k\in A}
\alpha_\xi^{m_k'-m_k}
\]
The set $\{\alpha_\xi^n:n\in\Z\}$ is discrete in $L^*$,
because $|\alpha_\xi|\ne 1$.
Therefore
\[
\alpha_\zeta^{qN}\frac{w_2}{w_1}=
\alpha_\xi^{m_k'-m_k}
\]
for all sufficiently large $k$ in $A$.

Recall that $q,N\ne 0$ and $w_1,w_2\in W$. It follows that
$\alpha_\zeta$ and $\alpha_\xi$ are multiplicatively dependent.
But we assumed the numbers
 $(\alpha_j)_{j\in S}$ to be multiplicatively
independent. Therefore $\xi=\zeta$.

By considering the quotient of the right hand side of eq.~(\ref{eq-tr-1}) and
its left hand side, we obtain:
\[
1= \lim_{k\to\infty,k\in A}\frac{\tilde w}{w_1a}\alpha_\xi^{m_k-qk}
\]
Therefore:
\begin{equation}\label{eq-3}
a=  \lim_{k\to\infty,k\in A}\frac{\tilde w}{w_1}\alpha_\xi^{m_k-qk}
\end{equation}
and consequently
\[
a= \frac{\tilde w}{w_1}\alpha_\xi^{m_k-qk}
\]
for all sufficiently large $k$ in $A$.

Together with eq.~(\ref{eq-tr-1}) this implies that 
\[
w_1\alpha_\xi^{qk} a 
= \tilde w\alpha_\xi^{m_k}
\text{ and }
 w_1^{-1}\alpha_\xi^{-qk}d 
= \tilde w^{-1}\alpha_\xi^{-m_k}
\]
Combining these two equalities we obtain $ad=1$.
Now recall that $\delta$ was an arbitrarily chosen element
in the intersection of $\Gamma_0$ with the Zariski open subset
\[
\Omega=\left\{
A=\begin{pmatrix}
a & c \\ b & d 
\end{pmatrix}
\in SL_2 : a,d\ne 0, A\gamma\ne\gamma A \right\}.
\]

Note that the condition $A\gamma=\gamma A$ implies that $A$ is a
diagonal matrix and therefore implies that $ad=1$.
 
Thus we have deduced:
{\em Every element of $\Gamma_0$ is contained
in the
algebraic subvariety 
\[
\left\{
\begin{pmatrix}
a & c \\ b & d 
\end{pmatrix}
\in SL_2 : ad=1 \text{ or }ad=0 \right\}.
\]}

But this contradicts the assumption that $\Gamma_0$ is Zariski-dense.
\end{proof}

\subsection{On the absolute values of eigenvalues}
For our main goal (i.e.~studying elliptic curves in quotients of
$SL_2(\C)$) we need only to consider the eigenvalues.
However, from the point of view of possible applications
to the study of geodesic length spectra of real hyperbolic
manifolds (see section 6 below) it might be interesting to deduce
a similar result for the absolute values of the eigenvalues.
This is the purpose of this subsection.
\begin{proposition}
Let $\Gamma$ be a subgroup of $SL_2(\C)$ which is dense
in the algebraic Zariski topology.

Then there exists infinitely many pairwise multiplicatively
independent positive real numbers 
which occur as the absolute value of an eigenvalue
for an element of $\Gamma$.
\end{proposition}
\begin{proof}
First we note that $|z|=\sqrt{z\bar z}$ for any complex number.
Using this fact, it is clear that for every finitely generated
subgroup $\Gamma$ of $SL_2(\C)$ there is a finitely
generated field $k$ such that every absolute value of
an eigenvalue for an element of $\Gamma$ is contained
in a finite extension field of degree at most $4$ over $k$:
We just have to take $k$ to be the extension field of $\Q$
generated by all the coefficients and their complex conjugates
for all elements in some fixed finite set of generators for $\Gamma$.

Thus the arguments in the proof of the preceding theorem
can be applied to deduce the following conclusion:

{\em
Either the statement of the proposition holds,

or (after conjugation with an appropriate element of $SL_2(\C)$)
we have
\[
\Gamma\subset \left\{ 
\begin{pmatrix}
a & b \\
c & d \\
\end{pmatrix}\in SL_2(\C) : |ad|\in\{0,1\}
\right\}.
\]
}

The condition $|ad|\in\{0,1\}$ is equivalent to 
$|ad|^2\in\{0,1\}$ which is a 
real algebraic condition.

Hence we have to discuss the {\em real} algebraic Zariski topology.
This is the topology whose closed sets are given as the zero
sets of polynomials in the complex coordinates {\em and} their
complex conjugates.

Since $\Gamma$ is Zariski-dense, the {\em real Zariski}-closure $S$
of $\Gamma$ in $SL_2(\C)$ is either the whole of $SL_2(\C)$ or
a real form of $SL_2(\C)$. 
Now $|ad|^2\in\{0,1\}$ defines a real
algebraic subset of $SL_2(\C)$.
Hence the real Zariski closure $S$ of $\Gamma$ cannot be the whole
of $SL_2(\C)$. Furthermore, since $\Gamma$ is discrete and infinite,
$S$ cannot be compact. Thus $S$ must be conjugate to $SL_2(\R)$.
However, this leads to a contradiction thanks to the
lemma below.
\end{proof}
\begin{lemma}
There is no element $A\in SL_2(\C)$ such that
\[
A\cdot SL_2(\R)\cdot A^{-1}
\subset 
 \left\{ 
\begin{pmatrix}
a & b \\
c & d \\
\end{pmatrix}\in SL_2(\C) : |ad|\in\{0,1\}
\right\}
\]
\end{lemma}
\begin{proof}
Let $\rho:SL_2(\C)\to\R$ denote the function given by
\[
\rho\begin{pmatrix}
a & b \\
c & d \\
\end{pmatrix} = |ad|.
\]
Now let us assume that the assertion of the lemma is wrong.
In other words: we assume that there exists an element

\[
A=
\begin{pmatrix}
x & y \\
z & w \\
\end{pmatrix}
\in SL_2(\C) 
\]
such that $\rho(g)\in\{0,1\}$ for every 
$g\in A\cdot SL_2(\R)\cdot A^{-1}$.

Since $SL_2(\R)$ is connected,
this implies that $\rho$ is constant, and its value either $0$ or $1$.

However, $\rho$ cannot be constantly zero, because $|ad|=0$
is equivalent to $ad=0$ and this is a complex algebraic condition.
Thus $\{g\in SL_2(\C):\rho(g)=0\}$ is a complex algebraic
subvariety and therefore cannot contain the group
$A\cdot SL_2(\R)\cdot A^{-1}$ which is dense in $SL_2(\C)$ with
respect to the complex Zariski topology.

This leaves the case where $\rho$ is constantly $1$.

Here explicit calculations show the following:
\[
\rho\left(
A \cdot
\begin{pmatrix}
1 & t \\
0 & 1 \\
\end{pmatrix}
\cdot A^{-1}
\right) = |1-(txz)^2|,
\]
and
\[
\rho\left(
A \cdot
\begin{pmatrix}
1 & 0 \\
t & 1 \\
\end{pmatrix}
\cdot A^{-1}
\right) = |1-(tyw)^2|
\]
Thus $xz$ and $yw$ are complex numbers with the property
that 
\[
|1-(txz)^2|=1=|1-(tyw)^2|
\] 
for {\em every} real number $t$.
This implies $xz=yw=0$.
But now
\begin{align*}
\rho\left(
A \cdot
\begin{pmatrix}
1 & 1 \\
-1 & 0 \\
\end{pmatrix}
\cdot A^{-1}
\right) &=|(xw-yw-xz)(-yz+yw+xz)| \\
&= |-xwyz| = 0 \ne 1 
\end{align*}
and we obtain a contradiction to the assumption 
that $\rho(ABA^{-1})=1$
for all $B\in SL_2(\R)$.
\end{proof}
\section{Equivalence of elliptic curves}
\subsection{Isogeny criteria}
An elliptic curve is a one-dimensional abelian variety,
or, equivalently a projective smooth algebraic curve of genus 1
(with a basepoint).
There are two natural equivalence relations between elliptic curves:
isomorphism (as algebraic variety) or isogeny.
Two varieties $V$ and $W$ are {\em isogenous} if there exists
a variety $Z$ and unramified coverings $\pi:Z\to V$, $\rho:Z\to W$.

Over the field of complex numbers, every elliptic curve can be realized
as the complex quotient manifold $\C/\left<1,\tau\right>_{\Z}$
where $\tau\in H^+=\{z:\Im(z)>0\}$.
Two elements $\tau,\tau'\in H^+$ define isomorphic resp.~isogenous
elliptic curves if both are contained  in the same $SL_2(\Z)$- resp.~%
$GL_2^+(\Q)$-orbit for the action on $H^+$ given by
\[
\begin{pmatrix}
a & b \\ c & d 
\end{pmatrix}
: z \mapsto \frac{az+b}{cz+d}.
\]
Here $GL_2^+(\Q)$ denotes the subgroup of $GL_2(\Q)$ containing
all elements with positive determinant.

We need some reformulations of these criteria.

\begin{lemma}\label{kern-crit}
Let $\Lambda$, $\Gamma$ be lattices in $\C$,
and $\Lambda_\Q=\Lambda\tensor\Q$, $\Gamma_\Q=\Gamma\tensor\Q$.
Consider the natural map
$\Phi:\Lambda_\Q\tensor_\Q\Gamma_\Q\to\C$
induced by the inclusion maps $\Gamma\hookrightarrow\C$, $\Lambda\hookrightarrow\C$.

Then $\C/\Lambda$ and $\C/\Gamma$ are isogenous iff $\dim_\Q\ker\Phi>0$.
\end{lemma}
\begin{proof}
We may assume $\Gamma=\left<1,\tau\right>_{\Z}$,
$\Lambda=\left<1,\sigma\right>_{\Z}$.
The kernel $\ker\Phi$ is positive-dimensional iff there
is a linear relation
\[
a + b \tau + c \sigma + d \tau\sigma =0
\]
with $(a,b,c,d)\in\Q^4\setminus\{(0,0,0,0)\}$.
Using $\sigma,\tau\in H^+$, one verifies that
\[
-\begin{pmatrix} a & b \\ c & d \end{pmatrix}
\in GL_2^+(\Q).
\]
Thus
\[
\sigma =-\frac{a+b\tau}{c+d\tau}
= -\begin{pmatrix} a & b \\ c & d \end{pmatrix}(\tau),
\]
i.e.\ $\dim\ker\Phi>0$ iff $\tau$ and $\sigma$ are contained
in the same $GL_2^+(\Q)$-orbit.
\end{proof}

\begin{lemma}\label{field-crit}
For a lattice $\left<\alpha,\beta\right>_{\Z}=
\Lambda\subset\C$ let $K_\Lambda$
denote the subfield of $\C$ given by
$K_\Lambda=\Q(\alpha/\beta)$.

Then $K_\Lambda$ depends only on $\Lambda$ and not of the choice
of the basis $(\alpha,\beta)$.

Let $\Lambda$ and $\tilde\Lambda$ be lattices in $\C$.

If $trdeg K_\Lambda/\Q>0$, then $\C/\Lambda$ and $\C/\tilde\Lambda $
are isogenous elliptic curves if and only if
$K_\Lambda=K_{\tilde \Lambda}$.
\end{lemma}
\begin{proof}
The independence of the choice of the basis is easily verified.

Furthermore, without loss of generality we may assume
$\Lambda=\left<1,\tau\right>$ and $\tilde\Lambda=\left<
1,\sigma\right>$ for some $\tau,\sigma\in H^+$.
Now the statement follows from the fact that
for transcendental complex numbers $\tau,\sigma$
we have $\Q(\tau)=\Q(\sigma)$ iff there are rational numbers
$a,b,c,d$ such that $\tau=(a+b\sigma)/(c+d\sigma)$.

Thus $\Q(\tau)=\Q(\sigma)$ iff $\sigma$ and $\tau$ are in the
same $GL_2^+(\Q)$-orbit in $H^+$.
\end{proof}

\subsection{Conjectures}
We now formulate a conjecture about an isogeny criterion for certain
elliptic curves:
\begin{conjecture}\label{our-conj}
Let $\alpha_1,\alpha_2\in\C$ be algebraic numbers with $|\alpha_i|>1$.
Let $E_i$ be the quotient manifold $\C^*/\{\alpha_i^k:k\in\Z\}$.

Then $E_1$ and $E_2$ are isogenous if and only if $\alpha_1$, $\alpha_2$
are multiplicatively dependent (in the sense of def.~\ref{def-mult-dep}).
\end{conjecture}

Note that $\C/\left<1,\tau\right>\simeq\C^*/\left<e^{2\pi i\tau}\right>$
for $\tau\in H^+$.
Thus for $\sigma,\tau\in H^+$ the condition ``$e^{2\pi i\sigma}$
and $e^{2\pi i\tau}$ are multiplicatively dependent'' translates into:
``There exists $a\in\Q^+$, $b\in\Q$ such that $\sigma=a\tau+b$''.

Therefore we can reformulate the above conjecture into terms of group
actions on the upper half plane $H^+$.

\begin{conjecture}
Let 
\[
B^+(\Q)=\left\{ \begin{pmatrix} a & b \\ 0 & a^{-1} \end{pmatrix}
 : a\in\Q^+, b\in \Q
\right \}
\]
and let $\sigma,\tau\in H^+$ be contained in the same $GL_2^+(\Q)$-orbit.

Assume that both $e^{2\pi i \sigma}$ and $e^{2\pi i\tau}$ are algebraic.
Then $\sigma$ and $\tau$ are already contained in the same $B^+(\Q)$-orbit.
\end{conjecture}

Next we prove that these two equivalent conjectures of ours are true
provided the famous Schanuels conjecture is right.

\begin{proposition}
Conjecture \ref{our-conj} holds,
if Schanuels conjecture is true.
\end{proposition}

Schanuels Conjecture is the 
far-reaching conjecture from
transcendental number theory
which encompasses many important conjecture in this area.
It states the following:

\begin{Sconjecture}
If $x_1,\ldots,x_n$ are $\Q$-linearly independent complex numbers,
then the transcendence degree of $\Q(x_1,\ldots,x_n,e^{x_1},\ldots,
e^{x_n})$ over $\Q$ is at least $n$.
\end{Sconjecture}

Now we prove the proposition.
\begin{proof}
Indeed, let $x_1=2\pi i$, $x_2=\log\alpha_1$, $x_3=\log\alpha_2$.
Schanuels conjecture then implies that either
\begin{enumerate}
\item
$\dim_\Q\left<2\pi i,\log\alpha_1,\log\alpha_2\right>\le 2$, 
or
\item
$2\pi i,\log\alpha_1,\log\alpha_2$ are all three algebraically independent.
\end{enumerate}

Since $\Re\log\alpha_i=\log|\alpha_i|>0$ (recall that we assumed
$|\alpha_i|>1$) for $i=1,2$, in the first case there exist
integers $n,m\in\Z\setminus\{0\}$ such that $\alpha_1^n=\alpha_2^m$,
i.e. $\alpha_1$ and $\alpha_2$ are multiplicatively dependent.

In the second case we can conclude
 that $\log\alpha_1/2\pi i$ and $\log\alpha_2/2\pi i$
are both transcendental and $\Q(\log\alpha_1/2\pi i)\ne
\Q(\log\alpha_2/2\pi i)$. Hence $\C^*/\left<\alpha_1\right>$ is not isogenous to
$\C^*/\left<\alpha_2\right>$ in this case.

Thus we have shown that either $\alpha_1$ and $\alpha_2$ are multiplicatively
dependent, or $\C^*/\left<\alpha_1\right>$ must be isogenous to
$\C^*/\left<\alpha_2\right>$.
\end{proof}
\begin{remark}
Actually we do not use Schanuels conjecture in its full strength,
but only a special case of it. However, even the special statement
we need is not yet proven.
\end{remark}

\section{Elliptic Curves in $SL_2(\C)/\Gamma$}
Let $\Gamma$ be a discrete cocompact subgroup of $SL_2(\C)$ and
$X=SL_2(\C)/\Gamma$ the quotient manifold.
We are interested in elliptic curves embedded into $X$.
Every elliptic curve embedded into $X$ is an orbit of
a reductive Lie subgroup $H$ of $SL_2(\C)$ with $H\simeq\C^*$
(see \cite{SMF}). Conversely, if $H$ is a Lie subgroup
of $SL_2(\C)$ with $H\simeq\C^*$ and 
$\#\left(H\cap\Gamma\right)=\infty$,
then $H/(H\cap\Gamma)$ is an elliptic curve 
embedded into $X$ as an $H$-orbit.
If $\gamma\in\Gamma$ is an element of infinite order
in a discrete cocompact subgroup $\Gamma$, then $\gamma$ is a semisimple
element of $SL_2(\C)$, and the connected component of the centralizer
\[
C(\gamma)=\{g\in SL_2(\C):g\gamma=\gamma g\}
\]
is such a Lie subgroup of $SL_2(\C)$ which has an elliptic curve
as a closed orbit in $X$.
Moreover this elliptic curve is isogenous to the quotient manifold
of $\C^*$ by the infinite cyclic subgroup generated by $\lambda$
where $\lambda$ is an eigenvalue of $\gamma\in SL_2(\C)$.

These facts (for which we refer to \cite{SMF})
establish the relationship between isogeny classes of
elliptic curves
embedded in $X$ on one side and eigenvalues of elements of $\Gamma$
on the other side.

\begin{proposition}
If conjecture \ref{our-conj} holds,
then for every discrete cocompact subgroup $\Gamma\subset SL_2(\C)$
there exist infinitely many isogeny classes of elliptic curves
embedded in $X=SL_2(\C)/\Gamma$.
\end{proposition}
\begin{proof}
If $\Gamma$ is discrete and cocompact in $SL_2(\C)$, then it must be
Zariski-dense.
Hence by thm.~\ref{thm-mult-ind} there are infinitely many
complex numbers $\lambda_1,\lambda_2,\ldots$ which are pairwise
multiplicatively independent
and which occur as eigenvalue for elements
$\gamma_1,\gamma_2,\ldots$ in $\Gamma$.

Being multiplicatively independent implies in particular that none
of these numbers $\lambda_i$ is a root of unity.

Furthermore, $\Gamma$ is conjugate to a subgroup of $SL_2(k)$
for some number field $k$ (see \cite{R}, Thm.~7.67), hence all the
numbers $\lambda_i$ are algebraic numbers.

Let $H_i$ be the centralizer of $\gamma_i$ in $SL_2(\C)$.
An element of $SL_2(\C)$ with an eigenvalue different from $1$ and $-1$
is semisimple. Hence $H_i\simeq\C^*$. 
Now $H_i\cap\Gamma$ is discrete and contains the element
$\gamma_i$. Because $\lambda_i$ is not a root of unity,
$\gamma_i$ is of infinite order.
It follows that $\left<\gamma_i\right>\simeq\Z$
and that $H_i/(\Gamma\cap H_i)$ is an elliptic curve which is isogenous
to $\C^*/\left<\lambda_i\right>$.

Thus the quotients $H_i/(\Gamma\cap H_i)$ are elliptic curves
embedded in $X=SL_2(\C)/\Gamma$ and, provided conj.~\ref{our-conj}
holds, these elliptic curves are pairwise non-isogenous since
the $\lambda_i$ are pairwise multiplicatively
independent.
\end{proof}

In particular:
\begin{corollary}\label{cor-schanuel}
If Schanuel's conjecture holds,
then for every discrete cocompact subgroup $\Gamma\subset SL_2(\C)$
there exists infinitely many isogeny classes of elliptic curves
embedded in $X=SL_2(\C)/\Gamma$.
\end{corollary}

\subsection{The case where $\Gamma\cap SL_2(\R)$ is Zariski dense}
\begin{theorem}\label{prop-real}
Let $\Gamma$ be a discrete subgroup of $SL_2(\C)$ and assume that
$\Gamma\cap SL_2(\R)$ is Zariski-dense in $SL_2$.

Then 
there exists infinitely many isogeny classes of elliptic curves
embedded in $X=SL_2(\C)/\Gamma$.
\end{theorem}
\begin{proof}
By thm.~\ref{thm-mult-ind} 
there are infinitely many pairwise multiplicatively
independent complex numbers $\lambda_i$
occuring as eigenvalues for
elements $\gamma\in\Gamma\cap SL_2(\R)$.
None of these $\lambda_i$ is a root of unity.

If $\lambda$ is an eigenvalue for a matrix $SL_2(\R)$, then
either $\lambda$ is real or $|\lambda|=1$.
If $\lambda$ is an eigenvalue for an element of a discrete
subgroup of $SL_2(\R)$ with $|\lambda|=1$, then $\lambda$
must be a root of unity.

Since none of the $\lambda_i$ is a root of unity, it follows
that all the numbers $\lambda_i$ are real.

Thus there are infinitely many elliptic curves $E_i$ in
$X=SL_2(\C)/\Gamma$ which are isogenous to $\C^*/\left<\lambda_i\right>$
where the numbers $\lambda_i$ are all real and pairwise
multiplicatively
independent.

We claim that at most two of these $E_i$ can be isogenous.
Assume the converse, i.e., let $\lambda_i$, $\lambda_j$ and $\lambda_k$
be pairwise multiplicatively independent real numbers larger than $1$
such that the three elliptic curves $E_i$, $E_j$ and $E_k$ are all isogenous.

Note that $E_i=\C/\left<2\pi i,\log\lambda_i\right>$
and similarily for $E_j$ and $E_k$.
Isogeny of $E_i$ and $E_j$ implies that there is a $\Q$-linear relation
between $4\pi^2$, $\log\lambda_i\log\lambda_j$, 
$2\pi i\log\lambda_i$ and $2\pi i\log\lambda_j$ (see lemma~\ref{kern-crit}).
Now $4\pi^2\in\R$ and ${\log\lambda_i\log\lambda_j\in\R}$,
while $2\pi i\log\lambda_i$ and $2\pi i\log\lambda_j$ are $\Q$-linearly independent
elements of $i\R$. Therefore a $\Q$-linear  relation can only exists
if $4\pi^2/(\log\lambda_i\log\lambda_j)\in\Q$.

Similarily the existence of an isogeny of between $E_j$ and $E_k$ implies
$4\pi^2/(\log\lambda_j\log\lambda_k)\in\Q$.

Combined, this yields $(\log\lambda_i\log\lambda_j)/(\log\lambda_j\log\lambda_k)=
\log\lambda_i/\log\lambda_k\in\Q$ which contradicts the assumption of
$\lambda_i$ and $\lambda_k$ being multiplicatively independent.

This proves the claim.

Thus we obtain 
an infinite family of elliptic curves in $SL_2(\C)/\Gamma$ such
that for each of these curves there is at most one other curve in this
family to which it is isogenous.
It follows that there are infinitely many
isogeny classes.
\end{proof}
\section{Existence of $\Gamma$ for which $\Gamma\cap SL_2(\R)$
is cocompact in $SL_2(\R)$}
From a differential geometric point of view the torsion-free
discrete cocompact subgroups of $SL_2(\C)$ are precisely
those groups which occur as fundamental group of compact real
hyperbolic threefolds $M$.
The condition that $\Gamma\cap SL_2(\R)$ is cocompact in $SL_2(\R)$
translates into the condition that there is a real hyperbolic
surface geodesically embedded into $M$.

However, we use a different point of view to show the existence
of such $\Gamma$. There is an arithmetic way to produce 
discrete cocompact subgroups in $SL_2(\C)$ which we employ.

This arithmetic construction (see e.g.\cite{V}) 
is the following:
Let $K$ be either $\Q$ or a totally imaginary quadratic extension
of $\Q$, $\bar K$ the unique archimedean completion of $K$,
$L/K$ a quadratic extension, $\lambda\in K^*$ such that
$\lambda\not\in N_{L/K}(L^*)$.
Then a central simple $K$-algebra can be defined by
$A=\{a+bt:a,b\in L\}$  with multiplication given
by $at=ta^{\sigma}$ (for $Gal(L/K)=\{id,\sigma\}$)
and $t^2=\lambda$. The elements of norm one
constitute a $K$-anisotropic simple $K$-group $S$.
Now $S({\mathcal O}_K)$ becomes a discrete cocompact
subgroup of $S(\bar K)$. If $\bar K=\R$, then $S(\bar K)=SL_2(\R)$
if $A\tensor\R\simeq Mat(2,\R)$ and $S(\bar K)=SU(2)$ if 
$A\tensor\R$ is isomorphic to the algebra of quaternions.

We use this in the following way:
Let $F_1=\Q[\sqrt{2}]$, $F_2=\Q[i]$, $F_3=\Q[i,\sqrt{2}]$
and $p=5$.
We observe that the prime ideal $(5)$ splits in $F_2$:
$5=(2+i)(2-i)$. Now $(2+i)$ is prime in $\Z[i]$
and both residue class fields for $5$ in $\Z$ resp.~$2+i$ (or $2-i$)
 in $\Z[i]$
are isomorphic to the finite field $\F_5=\Z/5\Z$. Note that $2$ is not
a square in $\F_5$.
As a consequence 
the prime ideals $(5)$ and $(2+i)$ (and similarily for $(2-i)$)
are totally inert
with respect to the
the field extensions
$\Q[\sqrt{2}]/\Q$ resp.\ $\Q[i,\sqrt{2}]/\Q[i]$.
It follows that $5$ is not contained in the image of the
{\em norm} for either the field extension $\Q[i,\sqrt{2}]/\Q[i]$
or the field extension $\Q[\sqrt{2}]/\Q$.
 
Thus we may use the above construction with
\[
(K,L,\lambda)=(\Q[i],\Q[i,\sqrt 2],5)
\] resp. $=(\Q,\Q[\sqrt 2],5)$ 
to obtain a discrete cocompact
subgroup $\Gamma$ resp. $\Gamma_1$ in $S(\C)\simeq SL_2(\C)$ resp. $S(\R)$.
Evidently $\Gamma_1=\Gamma\cap S(\R)$.
Now observe that $\Q[\sqrt{2}]\subset\R$ implies $A\tensor\R\simeq Mat(2,\R)$.
Thus $S(\R)\simeq SL_2(\R)$.

We have thus established:
\begin{proposition}
There exists a discrete subgroup $\Gamma$ in $SL_2(\C)$
such that both $SL_2(\C)/\Gamma$ and $SL_2(\R)/(SL_2(\R)\cap\Gamma)$
are compact.
\end{proposition}

In combination with thm.~\ref{prop-real} this implies
the following:
\begin{corollary}\label{cor-ex-real}
There exists a discrete cocompact subgroup $\Gamma$ in $SL_2(\C)$
such that the complex quotient manifold $X=SL_2(\C)/\Gamma$
contains infinitely many pairwise non-isogenous elliptic curves.
\end{corollary}

\section{Geodesic length spectra for hyperbolic manifolds}

Here we want to relate our results on eigenvalues of elements
of discrete groups to the study of closed geodesics 
on real hyperbolic manifolds (As standard references for
hyperbolic manifolds, see \cite{EGM},\cite{RH}).

A {\em real hyperbolic manifold} is a complete Riemannian manifold
with constant curvature $-1$. In each dimension $n$ there is a unique
simply-connected real hyperbolic manifold $H^n$.

Let $\H=\{z+wj:z,w\in\C\}$ denote the division algebra
of {\em quaternions}, i.e., the algebra given by $j^2=-1$
and $zj=j\bar z$ for all $z\in\C$.

Now $H^2$ can be realized as $H^2\simeq\{z+tj\in\H:z\in\R, t\in\R^+\}$
and $H^3$ as $H^3\simeq\{z+tj\in\H,z\in\C, t\in\R^+\}$.
In both cases the hyperbolic metric is obtained from the
euclidean metric by multiplication
with $1/t$.
Let $\rho$ denote the induced distance function.

The isometry group $G$ of $H^2$ resp. $H^3$ is $PSL_2(\R)$
resp. $PSL_2(\C)$ with the action given by
\[
\begin{pmatrix}
a & b \\ c & d 
\end{pmatrix}
:
\zeta \mapsto (a\zeta+b)(c\zeta+d)^{-1}
\]
where the calculations take place in the algebra of quaternions.

Explicit calculations show that for any $A\in G$ we have
\[
\inf_{x\in H}\rho(x,Ax)=\log(\max\{|\lambda|^2,|\lambda^{-2}|\})
\]
where the infimum is taken over all points of $H^2$ resp. $H^3$
and $(\lambda,\lambda^{-1})$ are the roots of the characteristic
polynomial of $\tilde A$ where $\tilde A$ is an element of $SL_2(\C)$
which projects onto $A\in G\subset PSL_2(\C)=SL_2(\C)/\{I,-I\}$.

For a complete Riemannian manifold with strictly negative curvature
there is a unique closed geodesic for every element of the
fundamental group.

Therefore:
If $\Gamma$ is a torsion-free discrete subgroup
of $G$ then the set of lengths of closed geodesics of $H/\Gamma$
coincides with the set of logarithms of absolute values of squares of
eigenvalues of elements of $\Gamma$.

Moreover, if $H=H^3$, one can show that the logarithm of the eigenvalue
of an element $g\in\Gamma$
is the ``complex length'' of the corresponding
closed geodesic in the following sense:
Let $\gamma$ be a closed geodesic in a compact hyperbolic $3$-fold $M$.
Let $s$ be the length of $\gamma$ in the usual sense.
If we fix a point $p\in\gamma$, then the holonomy along $\gamma$
defines an orthogonal transformation of the normal space
$T_p(M)/T_p(\gamma)$. This normal space is isomorphic to $\R^2$,
thus an orthogonal transformation is simply a rotation by an angle 
$\theta$. Now the ``complex length'' of $\gamma$ is defined to be
$s+i\theta$ (\cite{Rd}).

The set of all real resp.~complex numbers occuring as (complex) length for
a closed geodesic is denoted as (complex) geodesic length spectrum.
(In the literature, usually multiplicities are taken into account,
and sometimes only simple closed geodesics are considered.
For our point of interest (the $\Q$-linear independence of geodesic
lengths) these distinctions are of no relevance.)

Therefore we obtain:
\begin{proposition}
Assume that $M$ is a compact real hyperbolic $3$-manifold.
Then there exist infinitely many closed geodesics on $M$
such that their {\em complex lengths} are pairwise $\Q$-linearly independent.
\end{proposition}

Using the results of \S2.4. on the absolute values of the eigenvalues we also
obtain:

\begin{theorem}\label{thm-geodesic}
Let $M$ be a compact real hyperbolic manifold of dimension
two or three and
$\Lambda$ its geodesic length spectrum.

Then $\Lambda$ contains infinitely many pairwise
$\Q$-linearly independent elements.
\end{theorem}

Another consequence is the following:
\begin{corollary}
Let $\Gamma$ be a Zariski-dense subgroup in $SL_2(\C)$.

Then there exist two elements $\gamma_1,\gamma_2\in\Gamma$
with
respective eigenvalues $\lambda_1,\lambda_2\in\R$ such that
the numbers $|\log\lambda_1|,|\log\lambda_2|$ generate a dense
subgroup of the additive group $(\R,+)$.
\end{corollary}

There is a related  a result of Benoist (\cite{B}) which implies that
the subgroup of $(\R,+)$ generated by {\sl all} the logarithms
of the absolute values of eigenvalues of elements of $\Gamma$
is dense.
Thus, for $SL_2(\R)$ and $SL_2(\C)$ we can improve this result of Benoist.
However, Benoist's work applies to other semisimple Lie groups as well,
where our results concern only $SL_2(\R)$ and $SL_2(\C)$.

\end{document}